\documentclass[a4paper,12pt]{article}
\usepackage{enumerate}
\usepackage{latexsym}
\usepackage{amsfonts}
\usepackage[only,ninrm,elvrm,twlrm,sixrm,egtrm,tenrm]{rawfonts}
\usepackage{indentfirst}
\usepackage{amsmath}
\usepackage[noend]{algorithmic}
\usepackage{algorithm}
\usepackage{graphicx,psfrag}
\usepackage{graphics}
\usepackage{makeidx}
\newtheorem{thm}{Theorem}[section]
\newtheorem{cor}[thm]{Corollary}

\def\qed{\hfill \nopagebreak\rule{5pt}{8pt}}
\title{\bf Nonexistence of triples of nonisomorphic
connected graphs with isomorphic connected $P_3$-graphs \footnote{
Research supported by NSFC. } }
\author{
\small  Xueliang Li and Yan Liu\\
[2mm]
\small Center for Combinatorics and LPMC \\
\small Nankai University, Tianjin 300071, China \\
\small lxl@nankai.edu.cn\\ }
\date{ }
\begin{document}
\maketitle \vskip5mm
\begin{abstract}
In the paper "Broersma and Hoede, {\it Path graphs}, J. Graph
Theory {\bf 13} (1989) 427-444", the authors proposed a problem
whether there is a triple of mutually nonisomorphic connected
graphs which have an isomorphic connected $P_3$-graph. For a long
time, this problem remains unanswered. In this paper, we give it a
negative answer that there is no such triple, and thus completely
solve this problem.\\
[2mm] Keywords: path graph, connected, isomorphism
\end{abstract}

\section {Introduction}

Broersma and Hoede \cite{bh} generalized the concept of line
graphs to that of path graphs by defining adjacency as follows.
Let $k$ be a positive integer, and $P_k$ and $C_k$ denote a path
and a cycle with $k$ vertices, respectively. Let $\pi_k(G)$ be the
set of all $P_k$'s in $G$. The {\it path graph} $P_k(G)$ of $G$ is
a graph with vertex set $\pi_k(G)$ in which two $P_k$'s are
adjacent whenever their union is a path $P_{k+1}$ or a cycle
$C_k$. Broersma and Hoede got many results on $P_3$-graphs,
especially, described two infinite classes of pairs of
nonisomorphic connected graphs which have isomorphic connected
$P_3$-graphs. They also raised a number of unsolved problems or
questions, all of which have been solved during these year, but
only the following one remains unanswered.\\

\noindent{\bf Problem}. Whether there exists a triple of mutually
nonisomorphic connected graphs which have an isomorphic connected
$P_3$-graph ?\\

For $k=2$, i.e., line graphs, from Whitney's result (see
\cite{hb}) it is not difficult to see that the problem has a
negative answer. In \cite {lz} the authors showed that for $k\geq
4$ there are not only triples of but also arbitrarily many
mutually nonisomorphic connected graphs with isomorphic connected
$P_k$-graphs. However, interestingly we will show in this paper
that for $k=3$ there does not exist any triple of mutually
nonisomorphic connected graphs with an isomorphic connected
$P_3$-graph, just like the case for $k=2$ but very different from
the case for $k\geq 4$. Note that If one drops the connectedness
of the original graph or its $P_3$-graph, then it is easy to find
arbitrarily many mutually nonisomorphic graphs with an isomorphic
$P_3$-graph.

\section {Preliminaries}

All graphs in this paper are undirected, finite and simple. We
follow the terminology and notations used in \cite{aehj, bm}. If
$\sigma$ is an isomorphism from $G$ to $H$, then $\sigma$ induces
a $P_k$-isomorphism $\sigma^*$ from $G$ to $H$, where
$\sigma^*(a_1a_2\cdots a_k)=\sigma(a_1)\sigma(a_2)\cdots
\sigma(a_k)$ for all $a_1a_2\cdots a_k\in \pi_k(G)$. A
$P_k$-isomorphism $\tau$ is {\it induced} if $\tau=\sigma^*$ for
some isomorphism $\sigma$. If $\tau_i$ is a $P_k$-isomorphism from
$G_i$ to $H_i$ for $i=1$ and $2$, then we say that $\tau_1$ and
$\tau_2$ are {\it equivalent} if there are isomorphisms $\sigma$
and $\rho$ from $G_1$ to $G_2$ and $H_1$ to $H_2$, respectively,
such that $\tau_1=(\rho^*)^{-1}\circ \tau_2\circ \sigma^*$.

Define an {\it $i$-thorn} to be a $P_3$ with exactly $i$ ($i=1$ or
$2$) terminal ends in $G$. Let $T_i(G)$ be the set of $i$-thorns
in $G$. We say that two $P_3$-isomorphisms $\tau_i$ from $G_i$ to
$H_i$ for $i=1$ and $2$, are {\it $T$-related} if $(i)$ $G_1$ and
$G_2$ differ only in their star components, so do $H_1$ and $H_2$;
$(ii)$ $|T_2(G_1)|=|T_2(G_2)|$; and $(iii)$
$\tau_1(\alpha)=\tau_2(\alpha)$ for every $\alpha\in
\pi_3(G_1)-T_2(G_1)=\pi_3(G_2)-T_2(G_2)$.

Consider two $1$-thorns $abc$ and $abd$ where $deg(a)\geq 2$ and
$deg(c)=deg(d)=1$, then swapping $abc$ and $abd$ gives a
$P_3$-isomorphism, which we call a {\it $B$-swap}.

Suppose $abcde$ is a $P_5$ in $G$ such that both $abc$ and $cde$
are terminal $1$-thorns, i.e., $deg(a)=deg(e)=1$ and $deg(c)=2$,
then swapping $abc$ and $cde$ gives a $P_3$-isomorphism, which we
call an {\it $S$-swap}.

For distinct $a, b\in V(G)$, let $D_{a, b}$ denote the subgraph of
$G$ consisting of the union of all $P_3$'s with ends $a$ and $b$
and with middle vertex of degree $2$ in $G$. If $D_{a, b}$ is
nonempty we call it a {\it diamond} with ends $a$ and $b$. We
usually write $V(D_{a, b})-\{a, b\}$ as $\{c_1, c_2,\cdots, c_k\}$
and call $k$ the {\it width} of $D_{a,b}$, and refer to $D_{a, b}$
as a $k$-diamond. Note that if $a\sim b$, the edge $ab$ is not
included in $D_{a, b}$. To distinguish the two possibilities, we
say that the diamond $D_{a, b}$ is {\it braced} if $a\sim b$ and
{\it unbraced} otherwise. For $1\leq i<j\leq k$, the $P_3$'s
$ac_ib$ are called {\it diamond paths} while the pair of $P_3$'s
$c_iac_j$ and $c_ibc_j$ is called a {\it diamond pair}. Then
swapping $c_iac_j$ and $c_ibc_j$ gives a $P_3$-isomorphism, which
we call a {\it $D$-swap}.

Suppose $\tau_1$ and $\tau_2$ are $P_3$-isomorphisms from $G$ to
$H$. We say that $\tau_1$ and $\tau_2$ are {\it $B$-related} if
$\tau_2^{-1}\circ \tau_1$ is the identity or a composition of
$B$-swaps. The {\it $S$-related} and {\it $D$-related} are defined
similarly. We use joins of these four equivalence relations: for
example, two $P_3$-isomorphisms are {\it $TBSD$-related} if we can
get from one to the other by a chain of zero or more $T$-, $B$-,
$S$- and/or $D$-relations.

 The following is the main result of
\cite{aehj}, based on which we shall solve our problem by case
analysis.

\begin{thm}
Let $\tau$ be a $P_3$-isomorphism from $G$ to $H$ such that at
least one of $G$ or $H$ is connected. Then $\tau$ is one of the
following:

\begin{itemize}
\setlength{\itemsep}{-1.1mm}

\item[(i)] $T$-related to a $P_3$-isomorphism of generalized
$K_{3, 3}$ type;

\item[(ii)] of special Whitney type; \vskip-1mm

\item[(iii)] $D$-related to a $P_3$-isomorphism of Whitney type$
3, 4, 5$ or $6$;

\item[(iv)] $D$-related to a $P_3$-isomorphism of bipartite type;
or

\item[(v)] $TBSD$-related to an induced $P_3$-isomorphism.
\end{itemize}
The definition for each of the above types will be given in the
successive subsections.
\end{thm}

For solving our problem, in Theorem 2.1 we only need to consider
that the original graphs $G$ and $H$ are nonisomorphic connected
graphs with $T_2(G)=T_2(H)=\emptyset$. Below, we will analyze the
types in Theorem 2.1 case by case in details.

\subsection {Generalized $K_{3,3}$ type}

First, we introduce the following notation which is used in the
definition of generalized $K_{3,3}$ type. We write $(c, d)ab(e,
f)\mapsto uvwxu$ if $G$ contains the edges $ab$, $ac$, $ad$, $be$,
$bf$, $H$ contains the $C_4$ $uvwxu$, and $\tau$ maps $cab\mapsto
xuv$, $dab\mapsto vwx$, $abe\mapsto uvw$ and $abf\mapsto wxu$. We
also write $abc(d, e)\mapsto uvwxy$ if $G$ contains the edges $ab,
bc, cd, ce$, $H$ contains the $P_5$ $uvwxy$, and $\tau$ maps
$abc\mapsto vwx$, $bcd\mapsto uvw$ and $bce\mapsto wxy$. This
notation will be reversed (e.g., $abcda\mapsto (w,x)uv(y,z)$) as
needed. Then, define the generalized $K_{3, 3}$ type as follows:

Either $\tau$ or $\tau^{-1}$ as in the following cases (i) through
(vii), or any equivalent $P_3$-isomorphism, is said to be of {\it
generalized $K_{3,3}$ type}.

\begin{itemize}
\setlength{\itemsep}{-1.1mm}

\item[(i)] $(c, d)ab(e, f)\mapsto u_1v_1u_2v_2u_1$, and $cad$ and
$ebf$ map to $P_3$ components of $H$.

\item[(ii)] $(c, d)ab(e, f)\mapsto u_1v_1u_2v_2u_1$, $kebfh\mapsto
yv_3u_1(v_1, v_2)$, and $cad$ maps to a $P_3$ component.

\item[(iii)] $(c, d)ab(e, f)\mapsto u_1v_1u_2v_2u_1$, $(k, l)eb(a,
f)\mapsto u_1v_1u_2v_3u_1$, $(h, i)fb(a, e)\mapsto
u_1v_2u_2v_3u_1$, and $cad$, $kel$ and $hfi$ map to $P_3$
components.

\item[(iv)] $(c, d)ab(e, f)\mapsto u_1v_1u_2v_2u_1$, $ecadg\mapsto
xu_3v_1(u_1, u_2)$, and $cebfh\mapsto yv_3u_1(v_1, v_2)$. Note
that $G$ and $H$ are connected and isomorphic.

\item[(v)] $(c, d)ab(e, f)\mapsto u_1v_1u_2v_2u_1$, $ebfhe\mapsto
(v_1, v_2)u_1v_3(y, z)$, and $cad$ maps to $yv_3z$. Again $G$ and
$H$ are connected and isomorphic.

\item[(vi)] $(c, d)ab(e, f)\mapsto u_1v_1u_2v_2u_1$, $(c, d)eb(a,
f)\mapsto u_1v_1u_2v_3u_1$, $(h, i)fb(a, e)\mapsto
u_1v_2u_2v_3u_1$, $aceda\mapsto (w, x)u_3v_1(u_1,u_2)$, and $hfi$
maps to $wu_3x$. Again $G$ and $H$ are connected and isomorphic.

\item[(vii)] The construction on $K_{3,3}$; $G\cong H\cong
K_{3,3}$.
\end{itemize}

\noindent{\bf Remark 1}. For generalized $K_{3, 3}$ type, it is
easy to get the following results:
\begin{itemize}
\setlength{\itemsep}{-1.1mm}

\item[1.] For cases (i), (ii) and (iii), $G$ and $H$ are
nonisomorphic, but $H$ is not connected and there are isolated
vertices in $P_3(G)$ and $P_3(H)$.

\item[2.] For cases (iv) and (vii), $G$ and $H$ are connected with
$T_2(G)=T_2(H)=\emptyset$, but $G$ and $H$ are isomorphic.

\item[3.] For cases (v) and (vi), $G$ and $H$ are connected, but
are isomorphic and there are isolated vertices in $P_3(G)$ and
$P_3(H)$.
\end{itemize}
Thus there is no pair of nonisomorphic connected graphs with
isomorphic connected $P_3$-graphs in generalized $K_{3, 3}$ type.

\subsection {Special Whitney type}

Let $SW$ be the graph obtained by subdividing each edge of $K_{1,
3}$ exactly once, then $P_3(SW)\cong C_6$. $\tau$ is a
$P_3$-isomorphism from $SW$ to $C_6$, then we say $\tau$,
$\tau^{-1}$ or any equivalent $P_3$-isomorphism is of {\it special
Whitney type}.

It is clear that $SW$ and $C_6$ are two nonisomorphic connected
graphs with isomorphic connected $P_3$-graphs.

\subsection {Whitney type $3, 4, 5$ or $6$}

In this subsection, we begin with a general idea which will be
used here and in the next subsection. Suppose $F$ is a graph. A
{\it diamond inflation} of $F$ is a graph obtained by replacing
each edge $ab\in E(F)$ by an unbraced $s_{ab}$-diamond $D_{a, b}$
$(s_{a, b}\geq 1)$, and adding $t_a$ terminal edges incident with
each $a\in V(F)$ $(t_a\geq 0)$. Suppose $\varphi$ is an
edge-isomorphism between graphs $F$ and $F^{'}$, and suppose $I$
and $I^{'}$ are diamond inflations of $F$ and $F^{'}$,
respectively, with the following property: for every $ab\in E(F)$,
if $\varphi(ab)=uv$ then (i) $s_{uv}=s_{ab}$ and (ii)
$t_u+t_v=t_a+t_b$. Obtain $G$ and $H$ from $I$ and $I^{'}$,
respectively, by adding star components to one of them (if
necessary) to make the numbers of $2$-thorns equal. Then we can
define a $P_3$-isomorphism $\tau$ from $G$ to $H$ and say that
$\tau$ is a {\it diamond inflation} of $\varphi$.

\noindent {\bf Remark 2}. If $D_{a, b}$ is a nontrivial diamond
(i.e., $s_{a, b}>1$) in $G$, then there exists a unique and
nontrivial diamond $D_{u, v}$ in $H$ (see the proof in
\cite{aehj}).

The type in this subsection is related to Whitney's exceptional
edge-isomorphisms which is stated as follows:
\begin{thm}
[{\bf Whitney} \cite{hw}] Suppose that $\varphi$ is an
edge-isomorphism from $G$ to $H$ where $G$ and $H$ are both
connected. If $\varphi$ is not induced, then $i=|E(G)|=|E(H)|\in
\{3, 4, 5, 6\}$, $G$ and $H$ are isomorphic to $W_i$ and $W^{'}_i$
in some order, and $\varphi$ is equivalent to $\varphi_{i}$ or
$\varphi^{-1}_i$, where
\begin{itemize}
\setlength{\itemsep}{-1.1mm}

\item[(i)] $W_6\cong W^{'}_6\cong K_4$, with $V(W_6)= \{a, b, c,
d\}$, $V(W^{'}_6)= \{u, v, w, x\}$, and $\varphi_6$ maps
$ab\mapsto uv$, $ac\mapsto uw$, $ad\mapsto vw$, $bc\mapsto ux$,
$bd\mapsto vx$ and $cd\mapsto wx$;

\item[(ii)] $W_5=W_6-cd$, $W_5^{'}=W_6^{'}-wx$ and
$\varphi_5=\varphi_6|E(W_5)$;

\item[(iii)] $W_4=W_6-\{bd, cd\}$, $W_4^{'}=W_6^{'}-\{vx, wx\}$
and $\varphi_4=\varphi_6|E(W_4)$; and

\item[(iv)] $W_3=W_6-\{bc, bd, cd\}\cong K_{1, 3}$,
$W_3^{'}=W_6^{'}-x\cong K_3$, and $\varphi_3=\varphi_6|E(W_3)$.
\end{itemize}
\end{thm}
Then a $P_3$-isomorphism $\tau$ is said to be of {\it Whitney type
$i$} if $\tau$ or $\tau^{-1}$ is equivalent to a diamond inflation
of $\varphi_i$ as above for $i = 3, 4, 5, 6$.

Denote by $t_z$ the number of terminal edges incident with $z$ for
$z$ in $\{a, b, c, d\}$ or $\{u, v, w, x\}$. For Whitney type
$P_3$-isomorphisms, according to condition (ii) of Diamond
Inflation, gives one equation from each pair of corresponding
edges of the original Whitney graphs. Then there is a same
solution for all four types:
\begin{equation}
\begin{cases}

t_u=\frac 12(t_a+t_b+t_c-t_d)\\

t_v=\frac 12(t_a+t_b-t_c+t_d)\\

t_w=\frac 12(t_a-t_b+t_c+t_d)\\

t_x=\frac 12(-t_a+t_b+t_c+t_d) & \text{(except for type 3)}
\end{cases}
\end{equation}
Because we require connected $P_3$-graphs, in the above four
equations we must have $t_z=0$ or $1$ for every $z\in \{a, b, c,
d\}\cup \{ u, v, w, x\}$. We write $(t_a, t_b, t_c, t_d)\mapsto
(t_u, t_v, t_w, t_x)$. If $t_a$, $t_b$, $t_c$, $t_d= 0$ or $1$,
then we get the corresponding solutions for $t_u$, $t_v$, $t_w$,
$t_x$ by (1). For example: (1, 0, 0, 1) $\mapsto$ (0, 1, 1, 0)
denotes that $t_a=1$, $t_b=t_c=0$ and $t_d=1$ correspond to
solutions $t_u=0$, $t_v=t_w=1$ and $t_x=0$ by (1). So it is easy
to check that there are only the following eight cases satisfying
$t_z=0$ or $1$ for every $z\in \{a, b, c, d\}\cup \{ u, v, w,
x\}$:

\begin{itemize}

\setlength{\itemsep}{-1.1mm}

\item[(i)] (0, 0, 0, 0) $\mapsto$ (0, 0, 0, 0).

\item[(ii)] (1, 1, 1, 1) $\mapsto$ (1, 1, 1, 1) (except for type
$3$).

\item[(iii)] (1, 1, 0, 0) $\mapsto$ (1, 1, 0, 0) ($ab\mapsto uv$).

\item[(iv)] (1, 0, 1, 0) $\mapsto$ (1, 0, 1, 0) ($ac\mapsto uw$).

\item[(v)] (1, 0, 0, 1) $\mapsto$ (0, 1, 1, 0) ($ad\mapsto vw$).

\item[(vi)] (0, 1, 1, 0) $\mapsto$ (1, 0, 0, 1) ($bc\mapsto ux$)
(except for type $3$).

\item[(vii)] (0, 1, 0, 1) $\mapsto$ (0, 1, 0, 1) ($bd\mapsto vx$)
(except for type $3$ or $4$).

\item[(viii)] (0, 0, 1, 1) $\mapsto$ (0, 0, 1, 1) ($cd\mapsto wx$)
(except for type $3, 4$ or $5$).
\end{itemize}

If a $P_3$-isomorphism $\tau$ or $\tau^{-1}$ is equivalent to a
diamond inflation of $\varphi_i$ as above, and falls into one of
the above cases (i) through (viii), then $\tau$ is said to be of
{\it special Whitney type $i$} for $i= 3, 4, 5$ or $6$. Thus only
in special Whitney type $i$ for $i = 3, 4, 5$ or $6$, we can find
pairs of nonisomorphic connected graphs with isomorphic connected
$P_3$-graphs if we choose suitable diamond widths.

\subsection {Bipartite type}

First, we also introduce the definition of bipartite type. Start
with a positive integer $k$ and an arbitrary bipartite graph $F$
with at least one edge and with a bipartition $(A, B)$. Let $I$
and $I^{'}$ be different diamond inflations of $F$, where each
edge $e$ is inflated to a diamond of the same width $s_e$ both
times, but in producing $I$ each vertex $v$ has $t_v$ terminal
edges added, while in producing $I^{'}$ it has $t^{'}_v$ terminal
edges added. where
\begin{equation}  t_v^{'}=
\begin{cases}
     t_v-k & \text{if $v\in A$}\\
     t_v+k & \text{if $v\in B$}
\end{cases}
\end{equation}
Thus, we need $t_v\geq k$ for all $v\in A$. Let $\varphi$ be the
identity edge-isomorphism from $F$ to itself. Clearly $\varphi$,
$I$ and $I^{'}$ satisfy condition (i) of Diamond Inflation, and
condition (ii) is satisfied because each edge of $F$ has the form
$ab$ with $a\in A$ and $b\in B$, so that
$t_a^{'}+t_b^{'}=(t_a-k)+(t_b+k)=t_a+t_b$. We can therefore obtain
a $P_3$-isomorphism $\tau$ by diamond inflation; $\tau$ is in
general not induced. We say $\tau$ and $\tau^{-1}$, or any
equivalent $P_3$-isomorphisms, are of {\it bipartite type}.

This case is similar to the above Whitney type. Because we require
that the $P_3$-graphs of $I$ and $I^{'}$ are connected, we must
have $t_v$, $t_v^{'} = 0$ or $1$ for every $v\in A\cup B$. Since
$k\leq t_v (v\in A)$, we have $k=0$ or $1$. If $k=0$, then $I\cong
I^{'}$. If $k=1$, then $t_u=1$ for all $u\in A$ and $t_v=0$ for
all $v\in B$. Otherwise, if there is a vertex $u_{0}\in A$ with
$t_{u_0}=0$ or a vertex $v_{0}\in B$ with $t_{v_0}=1$, then
$t_{u_0}^{'}=-1$ or $t_{v_0}^{'}=2$ by (2). Therefore we have a
$P_3$-isomorphism $\tau_0$ from $I$ to $I^{'}$, where $t_u=1$ and
$t_u^{'}=0$ for all $u\in A$, $t_v=0$ and $t_v^{'}=1$ for all
$v\in B$, respectively. Then we say that $\tau_0$ and
$\tau_0^{-1}$, or any equivalent $P_3$-isomorphism, are of {\it
special bipartite type}. Therefore, this is the only case to find
pairs of nonisomorphic connected graphs which have isomorphic
connected $P_3$-graphs in the bipartite type.

\subsection {$TBSD$-related to an induced $P_3$-isomorphism}

In this subsection, we require that there is no isolated vertices
in $P_3$-graphs. Then all $P_3$-isomorphisms are $BSD$-related to
an induced one. It is clear that if two original graphs $G$ and
$H$ are connected with an isomorphic $P_3$-graph, then $G\cong H$
by the definition of $BSD$-related. Thus in this type, if we
require connected $P_3$-graphs, then the original graph and its
$P_3$-graph are one to one.

Then from the arguments in above five subsections, we can get the
following corollary which is essential to the solution of our
problem.
\begin{cor}
Let $\tau$ be a $P_3$-isomorphism from $G$ to $H$, where $G$ and
$H$ are nonisomorphic connected graphs with an isomorphic
connected $P_3$-graph. Then $\tau$ is one of the following:
\begin{itemize}
\setlength{\itemsep}{-1.1mm}

\item[(i)] of special Whitney type;

\item[(ii)] $D$-related to a $P_3$-isomorphism of special Whitney
type $3, 4, 5$ or $6$; or

\item[(iii)] $D$-related to a $P_3$-isomorphism of special
bipartite type.
\end{itemize}
\end{cor}

\section {Main result}

Now we can state and show the main result of this paper.
\begin{thm}
There is no triple of mutually nonisomorphic connected graphs with
an isomorphic connected $P_3$-graph.
\end{thm}
{\bf Proof.} Assume, to the contrary, that there exists a triple
of mutually nonisomorphic connected graphs $G_1$, $G_2$ and $G_3$
which have an isomorphic connected $P_3$-graph. Let $\tau_i$ be a
$P_3$-isomorphism from $G_i$ to $G_{i+1}$, then $\tau_i$ will be
one of three types in Corollary 2.3 for $i = 1, 2$.

{\bf Case 1.} $\tau_1$ and $\tau_2$ are of the same type.

{\bf Subcase 1.1} $\tau_1$ and $\tau_2$ are both of special
Whitney type.

Without loss of generality, let $G_1\cong SW$ and $G_2\cong C_6$.
Since $\tau_2$ is also of special Whitney type, it is clear that
$G_3\cong SW$. Thus $G_1\cong G_3$, a contradiction.

{\bf Subcase 1.2} $\tau_1$ and $\tau_2$ are both of $D$-related to
a $P_3$-isomorphism of special Whitney type $i$ for $i = 3, 4, 5$
or $6$.

Without loss of generality, we assume that $i=4$. Then $\tau_1$
and $\tau_2$ are $D$-related to a $P_3$-isomorphism of special
Whitney type $4$, and let $G_1$ and $G_2$ be diamond inflations of
$W_4$ and $W_4^{'}$, respectively, where $t_a=1$, $t_b=t_c=0$,
$t_d=1$, $t_u=0$, $t_v=t_w=1$ and $t_x=0$ (i.e., $ (1, 0, 0,
1)\mapsto (0, 1, 1, 0)$). Since $\tau_1$ and $\tau_2$ are of the
same type, $G_3$ is also a diamond inflation of $W_4$, and
$t_a=1$, $t_b=t_c=0$, $t_d=1$ by (1). Hence $G_1\cong G_3$, also a
contradiction.

{\bf Subcase 1.3} $\tau_1$ and $\tau_2$ are both of $D$-related to
a $P_3$-isomorphism of special bipartite type.

This subcase is similar to Subcase 1.2. Denote by $F$ an arbitrary
bipartite graph with a bipartition $(A, B)$. Then assume that
$G_1$ and $G_2$ are different diamond inflations of $F$,
respectively, where $t_u=1$ for all $u\in A$ and $t_v=0$ for all
$v\in B$ in $G_1$; $t_u=0$ for all $u\in A$ and $t_v=1$ for all
$v\in B$ in $G_2$. Thus we can easily obtain that $G_3$ is also a
diamond inflation of $F$ with $t_u=1$ for all $u\in A$ and $t_v=0$
for all $v\in B$ in $G_3$ by the definition of $\tau_2$. Then
$G_1\cong G_3$, contrary to the assumption.

{\bf Case 2.} $\tau_1$ and $\tau_2$ are of different types.

By the definition of special Whitney type, we know that it is also
a particular case of special Whitney type $3$, with the following
restrictions: (i) each edge $e$ with diamond width $s_e=1$ in
$K_{1, 3}$ and $K_3$, and (ii) $t_u=0$ for each vertex $u$ in
$K_{1, 3}$ and $K_3$. In fact, special Whitney type is the same as
special Whitney type $3$ in essence. Then in order to solve Case
2, we only need to distinguish the following two subcases:

{\bf Subcase 2.1} $\tau_1$ is of special Whitney type, and
$\tau_2$ is $D$-related to a $P_3$-isomorphism of special
bipartite type.

By the definition of $\tau_1$, $G_1$ or $G_2$ is a diamond
inflation of $W_3=K_{1, 3}$ or $W_3^{'}=K_3$; and also by
$\tau_2$, $G_2$ and $G_3$ are different diamond inflations of some
bipartite graph, respectively. Thus there is only one possibility:
$G_2$ is a diamond inflation of $K_{1, 3}$, where $K_{1, 3}$ has a
bipartition $A=\{a\}$, $B=\{b, c, d\}$. Then $G_1\cong C_6$ and
$G_2\cong SW$. It is easy to see that $t_a=t_b=t_c=t_d=0$ in
$K_{1, 3}$, a contradiction to the definition of special bipartite
type, where $t_a=1$ or $t_b=t_c=t_d=1$.

{\bf Subcase 2.2} $\tau_1$ is $D$-related to a $P_3$-isomorphism
of special Whitney type $i$ for $i = 3, 4, 5$ or $6$, and $\tau_2$
is $D$-related to a $P_3$-isomorphism of special bipartite type.

For $i = 4, 5$ or $6$, if $\tau_1$ is $D$-related to a
$P_3$-isomorphism of special Whitney type $i$, then $G_1$ and
$G_2$ are diamond inflations of $W_i$ and $W_i^{'}$ which have odd
cycles. $G_2$ and $G_3$ are different diamond inflations of some
bipartite graph by the definition of $\tau_2$. Then $\tau_1$ must
be $D$-related to a $P_3$-isomorphism of special Whitney type $3$.
By the same argument as in Subcase 2.1, we obtain that $G_2$ is a
diamond inflation of $K_{1, 3}$, where $K_{1, 3}$ has a
bipartition $A=\{a\}$, $B=\{b, c, d\}$. By the definition of
special Whitney type $3$, $\tau_1$ falls into one of the following
four cases: (0, 0, 0, 0) $\mapsto$ (0, 0, 0, 0), (1, 1, 0, 0)
$\mapsto$ (1, 1, 0, 0) ($ab\mapsto uv$), (1, 0, 1, 0) $\mapsto$
(1, 0, 1, 0) ($ac\mapsto uw$), or (1, 0, 0, 1) $\mapsto$ (0, 1, 1,
0) ($ad\mapsto vw$). However, by the definition of special
bipartite type, there are only two choices: either $t_a=0$,
$t_b=t_c=t_d=1$, or $t_a=1$, $t_b=t_c=t_d=0$. Finally, there does
not exist any graph $G_2$ that has common property of two
different types at the same time. So $\tau_1$ and $\tau_2$ must be
of the same type, a contradiction. The proof is thus complete.\qed

\end{document}